\documentclass{article}

\usepackage{amssymb}
\usepackage{amsmath}
\usepackage[dvips]{graphicx}

\begin{document}


\noindent \textbf{                                  }

\noindent \textbf{}

\noindent \textbf{}

\noindent \textbf{}
\begin{center}
\noindent \textbf{THE ONE-DIMENSHIONAL INVERSE WAVE }

\noindent \textbf{SPECTRAL PROBLEM WITH DISCONTINUOUS WAVE SPEED}
\end{center} \noindent \textbf{}
\begin{center}
\noindent \textbf{R.F. Efendiev }
\end{center}
\noindent \textbf{}
\begin{center}
\textbf{Institute Applied Mathematics, Baku State University,}

\textbf{Z.Khalilov, 23, AZ1148, Baku, Azerbaijan, }rakibaz@yahoo.com\textbf{ \textit{  }}
\end{center}
\textbf{\textit{}}

\textbf{\textit{}}

\indent{\underline{\textbf{ABSTRACT}}.}

\noindent

The inverse problem for the Sturm- Liouville operator with complex periodic potential and positive discontinuous coefficients on the axis is studied.  Main characteristics of the fundamental solutions are investigated, the spectrum of the operator is studied.  We give formulation of the inverse problem, prove a uniqueness theorem and provide a constructive procedure for the solution of the inverse problem.

\textbf{\underbar{Key words}}:  Discontinuous equation, Spectral singularities, Inverse spectral problem, Continuous spectrum.

\textbf{\underbar{MSC}}:   34A36, 34L05, 47A10, 47A70

\textbf{\underbar{INTRODUCTION.}}

\textbf{\underbar{}}

\noindent        We consider the differential equation

\begin{equation} \label{GrindEQ__1_}
-y''\left(x\right)+\, q\left(x\right)y\left(x\right)=\lambda ^{2} \rho \left(x\right)y\left(x\right)
\end{equation}
in the space $L_{2} \left(-\infty ,+\infty \right)$ where the prime denotes the derivative with respect to the space coordinate and assume that the potential  $q\left(x\right)$ is of the form

\begin{equation} \label{GrindEQ__2_}
q(x)=\sum _{n=1}^{\infty }q_{n} e^{inx}  ,
\end{equation}
the condition $\sum _{n=1}^{\infty }\left|q_{n} \right|^{2}
=q<\infty  $ is satisfies, $\lambda $ is a complex number, and

\begin{equation} \label{GrindEQ__3_}
\rho \left(x\right)=\left\{\, \, \begin{array}{c} {1} \\ {\beta ^{2} } \end{array}\right. \, \, \, \begin{array}{c} {for} \\ {for} \end{array}\, \, \, \, \begin{array}{c} {x\ge 0,} \\ {x<0,\, \beta \ne 1,\beta >0\, .} \end{array}
\end{equation}
This equation, in the frequency domain, describes the wave propagation in a nonhomogeneous medium, where $q(x)$- the restoring is force and ${\raise0.7ex\hbox{$ 1 $}\!\mathord{\left/{\vphantom{1 \rho \left(x\right)}}\right.\kern-\nulldelimiterspace}\!\lower0.7ex\hbox{$ \rho \left(x\right) $}} $ is the wave speed. The discontinuities in $\rho \left(x\right)$ correspond to abrupt changes in the properties of the medium in which the wave propagates.

\noindent If$\rho \left(x\right)=1$, then equation \eqref{GrindEQ__1_} is called the potential equation.  Especially, we would like to indicate that generalized Legendre  equation, degenerate   hyper-

\noindent geometrical equation, Bessel's equation and also Mathieu equation after suitable substitution coincide with potential equation \eqref{GrindEQ__1_}[1,p.374]-[2].

\noindent In regard to the problems with discontinuous coefficients, we remark that Sabatier and his co-workers [3-6] studied the scattering for the impedance-potential equation and the similar problems were intensively studied by many authors in different statements [7], [8], but for periodic complex potential they are considered for the first time.

\noindent Firstly potential \eqref{GrindEQ__2_} was considered by M.G.Gasymov[9]. Later in 1990 the results obtained in [8] were extended by Pastur L.A., Tkachenko V.A [10]. As a final remark we mention some related work of Guillemin , Uribe [11] and [12-14].

\noindent In this paper we will study the spectrum and also solve the inverse problem for singular non-self-adjoint operator.  As the coefficient allows bounded analytic continuation to the upper half-plane of the complex plane$z=x+it$, we can conduct detailed analysis of problem \eqref{GrindEQ__1_}-\eqref{GrindEQ__3_}.

\noindent The paper consists of three parts.

\noindent    In part 1 we study the properties of fundamental system of solutions of equation \eqref{GrindEQ__1_}. The spectrum of problem \eqref{GrindEQ__1_}-\eqref{GrindEQ__3_} is investigated in part 2. In part 3 we give a formulation of the inverse problem, prove a uniqueness theorem and provide a constructive procedure for the solution of the inverse problem

\noindent

\noindent

\textbf{\underbar{1.  REPRESENTATION OF FUNDAMENTAL SOLUTIONS.}}

\noindent \textbf{\underbar{}}

Here we study the solutions of the main equation

\noindent

\[-y''\left(x\right)+\, q\left(x\right)y\left(x\right)=\lambda ^{2} \rho \left(x\right)y\left(x\right)\]
that will be convenient in future.

\noindent  We first consider the solutions $f_{1}^{+} \left(x,\lambda \right)$ and $f_{2}^{+} \left(x,\lambda \right)$, determined by the conditions at infinity

\[\, \, \mathop{\lim }\limits_{Imx\to \infty } f_{1}^{+} \left(x,\lambda \right)e^{-i\lambda x} =1,\]

\[\, \, \, \mathop{\lim }\limits_{Imx\to \infty } f_{2}^{+} \left(x,\lambda \right)e^{i\beta \lambda x} =1.\]
We can prove the existence of these solutions if the condition
 $\sum_{n=1}^{\infty }\left|q_{n}\right|^{2}
  =\\=q<\infty$ is fulfilled for the potential. This will be unique restriction on
the potential and later on we'll consider it to be fulfilled.

\noindent \textbf{\underline{Theorem1}.  }Let \textit{$q(x)$} be
of the form \eqref{GrindEQ__2_} and  \textit{$\rho \left(x\right)$
}satisfy condition \eqref{GrindEQ__3_}. Then equation
\eqref{GrindEQ__1_} has special solutions of the form

\begin{equation} \label{GrindEQ__4_}
f_{1}^{+} (x,\lambda )=e^{i\lambda x} \left(1+\sum _{n=1}^{\infty }\frac{1}{n+2\lambda } \sum _{\alpha =n}^{\infty }V_{n\alpha } e^{i\alpha x}   \right)\, ,\, \, \, \, \, \, \, \, for\, \, x\ge 0,
\end{equation}

\begin{equation} \label{GrindEQ__5_}
f_{2}^{+} (x,\lambda )=e^{-i\lambda \beta x} \left(1+\sum _{n=1}^{\infty }\frac{1}{n-2\lambda \beta } \sum _{\alpha =n}^{\infty }V_{n\alpha } e^{i\alpha x}   \right)\, ,\, \, \, \, \, \, for\, \, \, x<0.
\end{equation}
where the numbers $V_{n\alpha } $ are determined from the following recurrent relations

\begin{equation} \label{GrindEQ__6_}
\alpha (\alpha -n)V_{n\alpha } +\sum _{s=n}^{\alpha -1}q_{\alpha -s} V_{ns} =0 ,\, \, \, \, \, \, \, \, \, \, \, 1\le n<\alpha ,
\end{equation}

\begin{equation} \label{GrindEQ__7_}
\alpha \sum _{n=1}^{\alpha }V_{n\alpha } +q_{\alpha } =0 ,
\end{equation}
 and the series

\begin{equation} \label{GrindEQ__8_}
\sum _{n=1}^{\infty }\frac{1}{n} \sum _{\alpha =n}^{\infty }\alpha \left|V_{n\alpha } \right|
\end{equation}
converges.

\noindent The proof of the theorem is similar to the proof of [10] and therefore we don't cite it here.

\noindent \textbf{\underline{Remark1}: } If $\lambda \ne
-\frac{n}{2} $, $\lambda \ne \frac{n}{2\beta } $ and $Im\lambda
\ge 0$, then$f_{1}^{+} \left(x,\lambda \right)\in L_{2}
\left(0,+\infty \right)$, $f_{2}^{+} \left(x,\lambda \right)\in
L_{2} \left(-\infty ,0\right)$.

\noindent Extending $f_{1}^{+} \left(x,\lambda \right)$and
$f_{2}^{+} \left(x,\lambda \right)$ as solutions of equation
\eqref{GrindEQ__1_} on \textit{$x<0$} and \textit{ $x\ge 0$
}respectively and using the conjunction conditions

\begin{equation} \label{GrindEQ__9_}
\begin{array}{l} {y\left(0+\right)=y(0-),} \\ {y'\left(0+\right)=y'(0-),} \end{array}
\end{equation}
we can prove the following lemma.

\textbf{\underline{Lemma 1:} } $f_{1}^{+} \left(x,\lambda \right)$
and $f_{2}^{+} \left(x,\lambda \right)$  may be extended as
solutions of equation \eqref{GrindEQ__1_} on \textit{$x<0$}
and\textit{$x\ge 0$, }respectively. Then we get

\[f_{2}^{+} \left(x,\lambda \right)=C_{11} \left(\lambda \right)f_{1}^{+} \left(x,\lambda \right)+C_{12} \left(\lambda \right)f_{1}^{-} \left(x,\lambda \right)\, \, \, for\, \, \, x\ge 0,\]

\[f_{1}^{+} \left(x,\lambda \right)=C_{22} \left(\lambda \right)f_{2}^{+} \left(x,\lambda \right)+C_{21} \left(\lambda \right)f_{2}^{-} \left(x,\lambda \right)\, ,\, \, \, \, \, \, \, for\, \, \, x<0\, \, ,\]

\noindent where

\noindent

\[f_{1,2}^{-} \left(x,\lambda \right)=f_{1,2}^{+} \left(x,-\lambda \right),\]

\begin{equation} \label{GrindEQ__10_}
C_{11} \left(\lambda \right)=\frac{W[f_{2}^{+} \left(0,\lambda \right),f_{1}^{-} \left(0,\lambda \right)]}{2i\lambda } ,
\end{equation}

\[C_{12} \left(\lambda \right)=\frac{W[f_{1}^{+} \left(0,\lambda \right),f_{2}^{+} \left(0,\lambda \right)]}{2i\lambda } ,\]

\begin{equation} \label{GrindEQ__11_}
C_{22} \left(\lambda \right)=\frac{1}{\beta } C_{11} \left(-\lambda \right),  C_{21} \left(\lambda \right)=-\frac{1}{\beta } C_{12} \left(\lambda \right).
\end{equation}
\textbf{\underline{Proof:}} It is easy to see that equation
\eqref{GrindEQ__1_} has fundamental solutions$f_{1}^{+}
\left(x,\lambda \right)$,$f_{1}^{-} \left(x,\lambda \right)$
($f_{2}^{+} \left(x,\lambda \right)$,$f_{2}^{-} \left(x,\lambda
\right)$)   for which

\[W\left[f_{1}^{+} (x,\lambda ),f_{1}^{-} (x,\lambda )\right]=2i\lambda ,\]

\[W\left[f_{2}^{+} (x,\lambda ),f_{2}^{-} (x,\lambda )\right]=2i\lambda \beta ,\]
 is satisfied\textbf{}

\noindent  Really, since $W[{\rm \; }f_{1}^{+} \left(x,\lambda \right),f_{1}^{-} \left(x,\lambda \right)]$ and $W[{\rm \; }f_{2}^{+} \left(x,\lambda \right),f_{2}^{-} \left(x,\lambda \right)]$ are independent of $x$ and the functions $f_{1}^{+} \left(x,\lambda \right)$,$f_{1}^{-} \left(x,\lambda \right)$ and $f_{2}^{+} \left(x,\lambda \right)$,$f_{2}^{-} \left(x,\lambda \right)$ allow holomorphic continuation on $x$ to upper and lower half-planes, respectively, the Wronskian coincides as $Imx\to \infty $. We can show that

\begin{equation} \label{GrindEQ__12_}
\mathop{\lim }\limits_{Imx\to \infty } f_{1}^{\pm \left(j\right)} \left(x,\lambda \right)e^{\mp i\lambda x} =\left(\pm i\lambda \right)^{j} \, \, \, \, \, \, \, \, j=0,1,
\end{equation}

\begin{equation} \label{GrindEQ__13_}
\mathop{\lim }\limits_{Imx\to \infty } f_{2}^{\pm \left(j\right)} \left(x,\lambda \right)e^{\mp \lambda x} =\left(\pm i\lambda \beta \right)^{j} \, \, \, \, \, \, \, \, j=0,1.
\end{equation}
So that

\[W\left[f_{1}^{+} (x,\lambda ),f_{1}^{-} (x,\lambda )\right]=2i\lambda ,\]

\[W\left[f_{2}^{+} (x,\lambda ),f_{2}^{-} (x,\lambda )\right]=2i\lambda \beta .\]
Then each solution of equation \eqref{GrindEQ__1_} may be represented as linear combinations of these solutions.

\[f_{2}^{+} \left(x,\lambda \right)=C_{11} \left(\lambda \right)f_{1}^{+} \left(x,\lambda \right)+C_{12} \left(\lambda \right)f_{1}^{-} \left(x,\lambda \right)\, \, \, for\, \, \, x\ge 0.\]

\[f_{1}^{+} \left(x,\lambda \right)=C_{22} \left(\lambda \right)f_{2}^{+} \left(x,\lambda \right)+C_{21} \left(\lambda \right)f_{2}^{-} \left(x,\lambda \right)\, ,\, \, \, \, \, \, \, for\, \, \, x<0\, \, ,\]
Using the conjunction conditions \eqref{GrindEQ__9_} it is easy to obtain the relation (10-11).

\noindent Let

\begin{equation} \label{GrindEQ__14_}
f_{n}^{\pm } (x)=\mathop{\lim }\limits_{\lambda \to \mp \frac{n}{2} } (n\pm 2\lambda )f_{1}^{\pm } (x,\lambda )=\sum _{\alpha =n}^{\infty }V_{n\alpha } e^{i\alpha x} e^{-i\frac{n}{2} x}  ,
\end{equation}
    It follows from relation \eqref{GrindEQ__6_} that if$V_{nn} \ne 0$, then $V_{n\alpha } \ne 0$ for all $\alpha >n$ and therefore$f_{n}^{\pm } (x)\ne 0$. Consequently, the points $\pm \frac{n}{2} ,\, n\in N$ are not singular points for$f_{1}^{\pm } \left(x,\lambda \right)$.

\noindent Then $W[f_{n}^{\pm } \left(x\right),f_{1}^{\mp } \left(x,\mp \frac{n}{2} \right)]=0$ and consequently the functions$f_{n}^{\pm } \left(x\right),f_{1}^{\mp } \left(x,\mp \frac{n}{2} \right)$, that are solutions of equation \eqref{GrindEQ__1_} for$\lambda =\pm \frac{n}{2} \, \, \, $, are linear dependent.

\noindent Therefore

\begin{equation} \label{GrindEQ__15_}
f_{n}^{\pm } \left(x\right)=V_{nn} f_{1}^{\mp } \left(x,\mp \frac{n}{2} \right),
\end{equation}

\noindent \textbf{\underbar{2.1.  SPECTRUM OF OPERATOR }}$L$.\textbf{\underbar{}}

Let $L$ be an operator generated by the operation  $\frac{1}{\rho \left(x\right)} \left\{-\frac{d^{2} }{dx^{2} } +q\left(x\right)\right\}$

\noindent in the space $L_{2} \left(-\infty ,+\infty ,\rho \left(x\right)\right)$.\textbf{}

\noindent To study the spectrums of the operator $L$ at first we calculate the kernel of the resolvent of the operator  \textbf{$\left(L-\lambda ^{2} I\right)$ }by means of general methods.

\noindent To construct the kernel of the resolvent of operator$L$, we consider the equation

\[-y''\left(x\right)+\, q\left(x\right)y\left(x\right)=\lambda ^{2} \rho \left(x\right)y\left(x\right)+f\left(x\right).\]

Here, $f\left(x\right)$ is an arbitrary function belonging to
 $L_{2} \left(-\infty ,+\infty \right)$.  Divide the plane $\lambda
$ into sectors

\[S_{k} =\{ k\pi <\arg \lambda <(k+1)\pi \} ,k=0,1.\]
When $\lambda \in S_{0} $ , we note that every solution of equation \eqref{GrindEQ__1_} can be written in the form

\begin{equation} \label{GrindEQ__16_}
y\left(x,\lambda \right)=C_{1} \left(x,\lambda \right)f_{1}^{+} \left(x,\lambda \right)+C_{2} \left(x,\lambda \right)f_{2}^{+} \left(x,\lambda \right)\, \, \, .
\end{equation}
Using the method of variation of constant, we obtain that

\[C_{1}^{'} \left(x,\lambda \right)=-\frac{1}{W[f_{1}^{+} ,f_{2}^{+} ]} \rho \left(x\right)f_{2}^{+} \left(x,\lambda \right)f\left(x\right)\]

\[C_{2}^{'} \left(x,\lambda \right)=\frac{1}{W[f_{1}^{+} ,f_{2}^{+} ]} \rho \left(x\right)f_{1}^{+} \left(x,\lambda \right)f\left(x\right)\]
By virtue of the condition  $y\left(x,\lambda \right)\in L_{2}
\left(-\infty ,+\infty \right)$, we find that

\[C_{2} \left(\infty ,\lambda \right)=C_{1} \left(-\infty ,\lambda \right)=0.     \]
Consequently, we have

\[C_{1} \left(x,\lambda \right)=\int _{-\infty }^{x}\frac{1}{W[f_{1}^{+} ,f_{2}^{+} ]} \rho \left(t\right)f_{2}^{+} \left(t,\lambda \right)f\left(t\right) dt\]

\[C_{2} \left(x,\lambda \right)=-\int _{x}^{\infty }\frac{1}{W[f_{1}^{+} ,f_{2}^{+} ]} \rho \left(t\right)f_{1}^{+} \left(t,\lambda \right)f\left(t\right) dt.\]
Substitute them in \eqref{GrindEQ__16_} we get

\[y\left(x,\lambda \right)=\int _{-\infty }^{\infty }R_{11} \left(x,t,\lambda \right) \rho \left(t\right)f\left(t\right)dt\]
where

\begin{equation} \label{GrindEQ__17_}
R_{11} \left(x,t,\lambda \right)=\frac{1}{W[f_{1}^{+} ,f_{2}^{+} ]} \left\{\begin{array}{c} {f_{1}^{+} \left(x,\lambda \right)f_{2}^{+} \left(t,\lambda \right)\, \, \, \, \, \, \, \, \, for\, \, t<x} \\ {f_{1}^{+} \left(t,\lambda \right)f_{2}^{+} \left(x,\lambda \right)\, \, \, \, \, \, \, \, \, for\, \, t>x} \end{array}\right.        \lambda \in S_{0} .
\end{equation}
Calculating analogously we can\textbf{ }construct the kernel of
the resolvent on the sector  $S_{1} $, namely

\begin{equation} \label{GrindEQ__18_}
R_{12} \left(x,t,\lambda \right)=\frac{1}{W[f_{1}^{-} ,f_{2}^{-} ]} \left\{\begin{array}{c} {f_{1}^{-} \left(x,\lambda \right)f_{2}^{-} \left(t,\lambda \right)\, \, \, \, \, \, \, \, \, for\, \, t<x} \\ {f_{1}^{-} \left(t,\lambda \right)f_{2}^{-} \left(x,\lambda \right)\, \, \, \, \, \, \, \, \, for\, \, t>x} \end{array}\right.        \lambda \in S_{1} .
\end{equation}

\noindent \textbf{ \underbar{Lemma 2.}  }The spectrum of the operator $L$ consist of continuous spectrum filling in the axis$ $$\{ 0\le \lambda <+\infty \} $\textit{ }on which there may exist spectral singularities coinciding with the numbers  $\frac{n}{2\beta } ,\, \, \, \frac{n}{2} ,\, \, \, n=1,2,3,...$

\noindent \textbf{\underbar{Proof}:  }

\noindent  It follows from (17-18) that the resolvent exists for
all complex values of $\lambda $ and may have poles on the real
axis. Let's  investigate these poles. First, we'll prove that the
operator $L$ has no negative eigenvalues. Assume opposite. Let
$k<0$ be eigenvalue of the operator $L$ with corresponding
eigenfunction.   Then from (17-18) it follows that the number
$\lambda _{0} =\left|k_{0} \right|^{{\raise0.7ex\hbox{$ 1
$}\!\mathord{\left/{\vphantom{1
2}}\right.\kern-\nulldelimiterspace}\!\lower0.7ex\hbox{$ 2 $}} }
\exp \left(i{\raise0.7ex\hbox{$ \pi
$}\!\mathord{\left/{\vphantom{\pi
2}}\right.\kern-\nulldelimiterspace}\!\lower0.7ex\hbox{$ 2 $}}
\right)=i\left|k_{0} \right|^{{\raise0.7ex\hbox{$ 1
$}\!\mathord{\left/{\vphantom{1
2}}\right.\kern-\nulldelimiterspace}\!\lower0.7ex\hbox{$ 2 $}} } $
should coincide with one of the numbers\textbf{ }$\frac{n}{2\beta
} ,\, \, \, \frac{n}{2} ,\, \, \, n=\pm 1,\pm 2,\pm 3,...$. The
obtained contradiction proves the absence of negative eigenvalues.

Now we'll investigate the function $R\left(x,t,\lambda \right)$ in
the neighborhood of $\lambda _{0} $ from$[0,\infty )$. Then the
number $\lambda _{0} $ coincides with one of  the
numbers$\frac{n}{2\beta } ,\, \, \, \frac{n}{2} ,\, \, \, n=\pm
1,\pm 2,\pm 3,...$. From (17-18) it follows that the limit
$\mathop{\lim \limits_{\lambda \to \lambda _{0} }\left(\lambda
-\lambda _{0} \right)} R\left(x,t,\lambda \right)=\\=R_{0}
\left(x,t\right)$ exists  and $R_{0} \left(x,t\right)$ is a
bounded function with respect to all the variables. Let $\theta
\left(x\right)$  be  an arbitrary finite function. Then $\varphi
\left(x\right)=\int _{-\infty }^{+\infty }R_{0} \left(x,t\right)
\theta \left(t\right)dt$ is a bounded solution of equation
\eqref{GrindEQ__1_} for$\lambda =\lambda _{0} $. Therefore
$\varphi \left(x\right)=C_{0} f_{1}^{+} \left(x,\lambda _{0}
\right)$. Comparison of the  last relation with formulae (17-18)
shows that if  $\lambda _{0} \ne \frac{n}{2} ,\, \, \lambda _{0}
\ne \frac{n}{2\beta } \, ,\, \, \, \, n\in N$ then   $C_{0} =0$and
so the kernel of  the resolvent has removable singularity at  the
point$\lambda _{0} $. So, it remains the $\lambda _{0} $ where has
poles of the first order. Since $f_{1}^{+} \left(x,\lambda _{0}
\right)\notin L_{2} \left(-\infty ,+\infty \right)$ then $\lambda
_{0}^{2} $ is a spectral singularity of the operator$L$.

Theorem is proved.

\textbf{\underbar{Corollary}:} The kernel $R\left(x,t,\lambda
\right)$ has no singularities on the axis$l_{1} =\{ \lambda :\arg
\lambda =\pi \} $, but on the axis $l_{0} =\{ \lambda :\lambda
>0\} $ it may have first order poles only at the points$\lambda
_{0} =\frac{n}{2} ,\, \, \lambda _{0} =\frac{n}{2\beta } \, ,\, \,
\, \, n\in N$. \noindent

\noindent\textbf {}

\noindent \textbf{\underbar{Lemma3:}} The coefficient $C_{12} \left(\lambda \right)$ is an analytic function in the $Im\lambda \ge 0$ and

\noindent has a finite number of simple zeros, moreover, if  $C_{12} \left(\lambda _{n} \right)=0$, then

\[\frac{d}{d\lambda } C_{12} \left(\lambda )\right|_{\lambda =\lambda _{n} }=-i\int _{-\infty }^{+\infty }\rho \left(x\right) f_{1}^{+} \left(x,\lambda _{n} \right)f_{2}^{+} \left(x,\lambda _{n} \right)dx.\]
\underline{\textbf{Proof}:}

\noindent From regularity $W[f_{1}^{+} \left(x,\lambda
\right),f_{2}^{+} \left(x,\lambda \right)]$ on $Im\lambda \ge 0$
and use the estimation

\[W[f_{1}^{+} \left(x,\lambda \right),f_{2}^{+} \left(x,\lambda \right)]=i\lambda \left(\beta -1\right)+O\left(\left|\lambda \right|^{-1} \right)\]
we get that  \textit{$C_{12} \left(\lambda \right)$ } has finite number of zeros.

\noindent Let's prove the second part of the theorem. By the standard method we can obtain from equation \eqref{GrindEQ__1_}

\noindent

\[f_{1}^{+} {''} \left(x,\lambda \right)\frac{d}{d\lambda } f_{2}^{+} \left(x,\lambda \right)-f_{1}^{+} \left(x,\lambda \right)\frac{d}{d\lambda } f_{2}^{+} {''}  \left(x,\lambda \right)=2\lambda \rho \left(x\right)f_{1}^{+} \left(x,\lambda \right)f_{2}^{+} \left(x,\lambda \right)\]

\noindent Integrating this equality from $-A$ to $x$ we find

\begin{equation} \label{GrindEQ__17_}
W[f_{1}^{+} \left(x,\lambda \right),\frac{d}{d\lambda } f_{2}^{+} \left(x,\lambda \right)]\left|_{-A}^{x} \right. =-2\lambda \int _{-A}^{x}\rho \left(x\right) f_{1}^{+} \left(x,\lambda \right)f_{2}^{+} \left(x,\lambda \right)dx.
\end{equation}
By the analogous way we get
\begin{equation} \label{GrindEQ__18_}
W[f_{1}^{+} \left(x,\lambda \right),\frac{d}{d\lambda } f_{2}^{+} \left(x,\lambda \right)]_{x}^{A} =-2\lambda \int _{x}^{A}\rho \left(x\right) f_{1}^{+} \left(x,\lambda \right)f_{2}^{+} \left(x,\lambda \right)dx.
\end{equation}
On the other hand, from Lemma1 follows that
\begin{equation} \label{GrindEQ__19_}
\frac{d}{d\lambda } (2i\lambda C_{12} \left(\lambda \right))=W[\frac{d}{d\lambda } f_{1}^{+} \left(x,\lambda \right),f_{2}^{+} \left(x,\lambda \right)]+W[f_{1} \left(x,\lambda \right),\frac{d}{d\lambda } f_{2}^{+} \left(x,\lambda \right)]
\end{equation}
\noindent Let's $\lambda =\lambda _{n} $ is one of the zeros of
the $C_{12} \left(\lambda \right)$. Comparing the formulas (19-21)
we get
\begin{equation} \label{GrindEQ__20_}
\begin{array}{l} {\left(2i\lambda \frac{d}{d\lambda } C_{12} \left(\lambda \right)+2iC_{12} \left(\lambda \right)\right)\left|_{\lambda =\lambda _{n} } \right. =2\lambda _{n} \int _{-A}^{A}\rho \left(x\right) f_{1}^{+} \left(x,\lambda _{n} \right)f_{2}^{+} \left(x,\lambda _{n} \right)dx+} \\ {+W[\frac{d}{d\lambda } f_{1}^{+} \left(x,\lambda _{n} \right),f_{2}^{+} \left(x,\lambda _{n} \right)]\left|_{x=-A} \right. +W[f_{1}^{+} \left(x,\lambda _{n} \right),\frac{d}{d\lambda } f_{2}^{+} \left(x,\lambda _{n} \right)]\left|_{x=A} \right. } \end{array}
\end{equation}
\noindent As the functions $f_{1}^{+} \left(x,\lambda _{n}
\right)$ and $f_{2}^{+} \left(x,\lambda _{n} \right)$ belong to
$L_{2} \left(-\infty ,+\infty \right)$ at $\lambda =\lambda _{n} $
therefore the  second and the third addends at the right hand side
in \eqref{GrindEQ__20_} are equal to zero at $A\to +\infty $, we
find that

\[\left. \frac{d}{d\lambda } C_{12} \left(\lambda \right)\right|_{\lambda =\lambda _{n} } =-i\int _{-\infty }^{+\infty }\rho \left(x\right) f_{1} \left(x,\lambda _{n} \right)\varphi _{2} \left(x,\lambda _{n} \right)dx.\]
Lemma3 is proved\textbf{}

\noindent \textbf{ }  For solutions  $f_{1}^{\pm } \left(x,\lambda \right)$ and $f_{2}^{\pm } \left(x,\lambda \right)$  we can obtain the asymptotic equalities

\noindent $ $$f_{1}^{\pm \left(j\right)} \left(0,\lambda \right)=\pm \left(i\lambda \right)^{j} +o\eqref{GrindEQ__1_}$             for$\left|\lambda \right|\to \infty ,\, \, \, j=0,1$,

\noindent

\noindent                              $f_{2}^{\pm \left(j\right)} \left(0,\lambda \right)=\pm \left(i\lambda \beta \right)^{j} +o\left(1\right)$          for $\left|\lambda \right|\to \infty ,\, \, \, j=0,1$

\noindent

\noindent  For simplicity we prove the first equality.

\noindent Since

\[f_{1}^{\pm } \left(0,\lambda \right)=1+\sum _{n=1}^{\infty }\sum _{\alpha =n}^{\infty }\frac{V_{n\alpha } }{n\pm 2\lambda }   \]
that

\[\left|f_{1}^{\pm } \left(0,\lambda \right)\right|\le 1+\sum _{n=1}^{\infty }\sum _{\alpha =n}^{\infty }\frac{\left|V_{n\alpha } \right|}{\left|n+2\lambda \right|} \le   \, 1+\sum _{n=1}^{\infty }\sum _{\alpha =n}^{\infty }\frac{\left|V_{n\alpha } \right|}{\sqrt{\left(n+2Re\lambda \right)^{2} +4Im^{2} \lambda } } \le   \, \, 1+\frac{1}{\left|Im\lambda \right|} \sum _{n=1}^{\infty }\sum _{\alpha =n}^{\infty }\frac{\alpha \left|V_{n\alpha } \right|}{n}   .\]
Therefore, as $\left|\lambda \right|\to \infty $, we obtain $f_{1}^{\pm } \left(0,\lambda \right)=1+o\left(1\right)$.

\noindent Analogously we can prove the rest asymptotic equalities as$\left|\lambda \right|\to \infty $,

\noindent Then for the coefficients $C_{12} (\lambda ),\, \, \, C_{12} (-\lambda ),\, $ we get the following asymptotic equalities

\begin{equation} \label{GrindEQ__21_}
C_{12} \left(\lambda \right)=\frac{1}{2i\lambda } \left(i\lambda \beta -i\lambda \right)+o\left(1\right)=-\frac{\beta +1}{2} +o\left(1\right),
\end{equation}

\[C_{12} \left(-\lambda \right)=-\frac{\beta +1}{2} +o\left(1\right),\]
These asymptotic equalities and analytical properties of the coefficients $C_{12} (\lambda ),\, \, \, C_{12} (-\lambda )$ make valid the following statement.

\noindent \textbf{\underline{Lemma 4}.}  The eigenvalues of
operator $L$are finite and coincide with zeros of the functions
$C_{12} (\lambda ),\, \, \, C_{12} (-\lambda )\, $ from sectors

\[S_{k} =\{ k\pi <\arg \lambda <(k+1)\pi \} ,k=0,1\]
respectively.

\noindent \textbf{\underline{Remark:}}  {Take into account
\eqref{GrindEQ__21_} we can obtain the useful on later relation}

\begin{equation} \label{GrindEQ__22_}
\beta =-2\mathop{\lim }\limits_{Im\lambda \to \infty } C_{12} \left(\lambda \right)-1
\end{equation}

\noindent \textbf{\underline{2.2. EIGENFUNCTION EXPANSIONS}.}

\noindent \textbf{\underline{Definition 2}}.  The points at which
resolvent have poles are called the singular numbers of
operator$L$.

\noindent Let  $\lambda _{1} ,\lambda _{2} ,....\lambda _{l} ,\lambda _{l+1} .....\lambda _{n} ...$be the singular numbers of operator $L$.At that

\[Re\lambda _{j} Im\lambda _{j} \ne 0,\, \, \, \, \, \, \, j=1,2,....l\]

\[Re\lambda _{j} Im\lambda _{j} =0,\, \, \, \, \, \, \, j=l+1,....n,...\]
The numbers $\lambda _{j} ,\, \, \, j=l+1,....n,..$are called the
spectral singularities of operator$L$. From the form of resolvent
it is easy to see that it has singular numbers (i.e. eigenvalues)
$\lambda _{1,} \lambda _{2,} ....\lambda _{l,} $in zeros of the
functions  $C_{12} \left(\lambda \right),\, \, C_{12}
\left(-\lambda \right)$ in the sectors $S_{k} =\{ k\pi <\arg
\lambda <(k+1)\pi \} ,k=0,1$  respectively. It directly follows
from Lemma 2 and representation (17-18) that kernel of resolvent
may have spectral singularities coinciding with the
numbers$\frac{n}{2\beta } ,\, \, \, \frac{n}{2} ,\, \, \,
n=1,2,3,...$.  Consequently taking \eqref{GrindEQ__15_} into
account we calculate
\begin{equation} \label{GrindEQ__23_}
\begin{array}{l} {\mathop{\lim }\limits_{\lambda \to {\raise0.7ex\hbox{$ n $}\!\mathord{\left/{\vphantom{n 2}}\right.\kern-\nulldelimiterspace}\!\lower0.7ex\hbox{$ 2 $}} } \left(n-2\lambda \right)R_{11} \left(x,t,\lambda \right)=\mathop{\lim }\limits_{\lambda \to {\raise0.7ex\hbox{$ n $}\!\mathord{\left/{\vphantom{n 2}}\right.\kern-\nulldelimiterspace}\!\lower0.7ex\hbox{$ 2 $}} } \left(n-2\lambda \right)\frac{1}{2i\lambda } [f_{1}^{+} \left(x,\lambda \right)f_{1}^{+} \left(t,\lambda \right)\frac{W[f_{2}^{+} ,f_{1}^{-} ]}{W[f_{1}^{+} ,f_{2}^{+} ]} +} \\ {+f_{1}^{+} \left(x,\lambda \right)f_{1}^{-} \left(t,\lambda \right)]=\frac{1}{in} [V_{nn} f_{1}^{+} \left(x,\frac{n}{2} \right)f_{1}^{+} \left(t,\frac{n}{2} \right)+} \\ {+V_{nn} f_{1}^{+} \left(x,\frac{n}{2} \right)f_{1}^{+} \left(t,\frac{n}{2} \right)]=\frac{2}{in} V_{nn} f_{1}^{+} \left(x,\frac{n}{2} \right)f_{1}^{+} \left(t,\frac{n}{2} \right).} \end{array}
\end{equation}
Analogously taking into account the denotation$\tilde{f}_{2}^{+}
\left(x,\lambda \right)=f_{2}^{+} \left(x,\lambda
\right)\left(n-2\lambda \beta \right)$, therewith, the function
$\tilde{f}_{2}^{+} \left(x,\lambda \right)$ has no poles at the
points$\lambda =\frac{n}{2\beta } ,\, \, \, n\in N$, \noindent we
get
\[
\begin{array}{l}
 \mathop {\lim }\limits_{\lambda  \to {\raise0.7ex\hbox{$n$} \!\mathord{\left/
 {\vphantom {n {2\beta }}}\right.\kern-\nulldelimiterspace}
\!\lower0.7ex\hbox{${2\beta }$}}} \left( {n - 2\lambda \beta }
\right)R_{11} \left( {x,t,\lambda } \right) = \mathop {\lim
}\limits_{\lambda  \to {\raise0.7ex\hbox{$n$} \!\mathord{\left/
 {\vphantom {n {2\beta }}}\right.\kern-\nulldelimiterspace}
\!\lower0.7ex\hbox{${2\beta }$}}} \left( {n - 2\lambda \beta } \right)\frac{1}{{W\left[ {f_1^ +  ,f_2^ +  } \right]}} \times  \\
  \times \left[ {C_{22} f_2^ +  \left( {x,\lambda } \right) + C_{21} \left( \lambda  \right)f_2^ -  \left( {x,\lambda } \right)} \right]f_2^ +  \left( {t,\lambda } \right) = \mathop {\lim }\limits_{\lambda  \to {\raise0.7ex\hbox{$n$} \!\mathord{\left/
 {\vphantom {n {2\beta }}}\right.\kern-\nulldelimiterspace}
\!\lower0.7ex\hbox{${2\beta }$}}} \left( {n - 2\lambda \beta } \right)\frac{1}{{W\left[ {f_1^ +  ,f_2^ +  } \right]}} \times  \\
  \times \left[ { - \frac{1}{\beta }\frac{{W\left[ {f_2^ -  ,f_1^ +  } \right]}}{{2i\lambda }}f_2^ +  \left( {x,\lambda } \right) + \frac{1}{\beta }\frac{{W\left[ {f_1^ +  ,f_2^ +  } \right]}}{{2i\lambda }}f_2^ -  \left( {x,\lambda } \right)} \right]f_2^ +  \left( {t,\lambda } \right) =  \\
  = \mathop {\lim }\limits_{\lambda  \to {\raise0.7ex\hbox{$n$} \!\mathord{\left/
 {\vphantom {n {2\beta }}}\right.\kern-\nulldelimiterspace}
\!\lower0.7ex\hbox{${2\beta }$}}} \left( {n - 2\lambda \beta } \right)\frac{{n - 2\lambda \beta }}{{W\left[ {f_1^ +  ,\tilde f_2^ +  } \right]}}\left[ { - \frac{1}{\beta }\frac{{W\left[ {f_2^ -  ,f_1^ +  } \right]}}{{2i\lambda }}\frac{{\tilde f_2^ +  \left( {x,\lambda } \right)}}{{n - 2\lambda \beta }} + } \right. \\
 \left. { + \frac{1}{\beta }\frac{{W\left[ {f_1^ +  ,\tilde f_2^ +  } \right]}}{{2i\lambda }}\frac{1}{{n - 2\lambda \beta }}f_2^ -  \left( {x,\lambda } \right)} \right]\frac{{\tilde f_2^ +  \left( {t,\lambda } \right)}}{{n - 2\lambda \beta }} = [ - \frac{1}{{in}}\frac{{W\left[ {f_2^ -  \left( {0,{\raise0.7ex\hbox{$n$} \!\mathord{\left/
 {\vphantom {n {2\beta }}}\right.\kern-\nulldelimiterspace}
\!\lower0.7ex\hbox{${2\beta }$}}} \right),f_1^ +  \left(
{0,{\raise0.7ex\hbox{$n$} \!\mathord{\left/
 {\vphantom {n {2\beta }}}\right.\kern-\nulldelimiterspace}
\!\lower0.7ex\hbox{${2\beta }$}}} \right)} \right]}}{{W\left[
{f_1^ +  \left( {0,{\raise0.7ex\hbox{$n$} \!\mathord{\left/
 {\vphantom {n {2\beta }}}\right.\kern-\nulldelimiterspace}
\!\lower0.7ex\hbox{${2\beta }$}}} \right),\tilde f_2^ +  \left(
{0,{\raise0.7ex\hbox{$n$} \!\mathord{\left/
 {\vphantom {n {2\beta }}}\right.\kern-\nulldelimiterspace}
\!\lower0.7ex\hbox{${2\beta }$}}} \right)} \right]}}\tilde f_2^ +
\left( {x,{\raise0.7ex\hbox{$n$} \!\mathord{\left/
 {\vphantom {n {2\beta }}}\right.\kern-\nulldelimiterspace}
\!\lower0.7ex\hbox{${2\beta }$}}} \right) +  \\
  + \frac{{W\left[ {f_1^ +  \left( {0,{\raise0.7ex\hbox{$n$} \!\mathord{\left/
 {\vphantom {n {2\beta }}}\right.\kern-\nulldelimiterspace}
\!\lower0.7ex\hbox{${2\beta }$}}} \right),\tilde f_2^ +  \left(
{0,{\raise0.7ex\hbox{$n$} \!\mathord{\left/
 {\vphantom {n {2\beta }}}\right.\kern-\nulldelimiterspace}
\!\lower0.7ex\hbox{${2\beta }$}}} \right)} \right]}}{{in}}f_2^ -
\left( {x,{\raise0.7ex\hbox{$n$} \!\mathord{\left/
 {\vphantom {n {2\beta }}}\right.\kern-\nulldelimiterspace}
\!\lower0.7ex\hbox{${2\beta }$}}} \right)]\tilde f_2^ +  \left(
{t,{\raise0.7ex\hbox{$n$} \!\mathord{\left/
 {\vphantom {n {2\beta }}}\right.\kern-\nulldelimiterspace}
\!\lower0.7ex\hbox{${2\beta }$}}} \right) = F\left( {x,t} \right) \\
 \end{array}
\]
\noindent \textbf{\underline{Lemma 6}:}   Let  $\psi
\left(x\right)$ be an arbitrary twice continuously differentiable
function belonging to$L_{2} \left(-\infty ,+\infty ,\rho
\left(x\right)\right)$ . Then
\[\int _{-\infty }^{+\infty }R\left(x,t,\lambda \right) \rho \left(t\right)\psi \left(t\right)dt=-\frac{\psi \left(x\right)}{\lambda ^{2} } +\frac{1}{\lambda ^{2} } \int _{-\infty }^{+\infty }R\left(x,t,\lambda \right) g\left(t\right)dt,\]
  where
\[g\left(t\right)=-\psi ''\left(x\right)+q\left(x\right)\psi \left(x\right)\in L_{2} \left(-\infty ,+\infty \right).\]
Integrating the both hand side along the circle $\left|\lambda \right|=R$ and passing to limit as $R\to \infty $   we get
\[\psi \left(x\right)=-\mathop{\lim }\limits_{R\to \infty } \frac{1}{2\pi i} \oint _{\left|\lambda \right|=R}2\lambda d\lambda \int _{-\infty }^{+\infty }R\left(x,t,\lambda \right)\rho \left(t\right)\psi \left(t\right)dt  \]
The function  $\int _{-\infty }^{+\infty }R\left(x,t,\lambda
\right)\rho \left(t\right)\psi \left(t\right)dt $ is analytical
inside the contour, with respect to $\lambda $ excepting the
points$\lambda =\lambda _{n} ,n=1,2,...l,\, \, \lambda
=\frac{n}{2} ,\, \, \lambda =\frac{n}{2\beta } ,\, \, \,
n=1,2,...$.   Denote by  $\Gamma_{0}^{+} \left(\Gamma_{0}^{-}
\right)$ the contour formed by segments  $[0,\frac{1}{2\beta }
-\delta ],[\frac{1}{2\beta } +\delta ,\frac{1}{2} -\delta
],...[\frac{n}{2\beta } +\delta ,\frac{n}{2} -\delta ]$   and
semi-circles of radius $\delta $ with the centers at points
$\frac{n}{2} ,\, \frac{n}{2\beta } \, \, \, n=1,2,...$ located in
upper (lower) half plane.

\noindent Then

\[\begin{array}{l} {\psi \left(x\right)=-\frac{1}{2i\pi } \int _{-\infty }^{+\infty }2\lambda \rho \left(t\right)\psi \left(t\right) [\int _{\Gamma_{0}^{+} }R_{11} \left(x,t,\lambda \right) d\lambda -\int _{\Gamma_{0}^{-} }R_{12} \left(x,t,\lambda \right) d\lambda ]dt=-\frac{1}{2i\pi } \int _{-\infty }^{+\infty }2\lambda \rho \left(t\right)\psi \left(t\right) \times } \\ {\times \int _{\Gamma_{0}^{-} }[R_{11} \left(x,t,\lambda \right) -R_{12} \left(x,t,\lambda \right)]d\lambda dt\mathop{+Res}\limits_{\lambda =\lambda _{n} } R_{11} \left(x,t,\lambda \right)\mathop{+Res}\limits_{\lambda =\frac{n}{2\beta } } R_{11} \left(x,t,\lambda \right)\mathop{+Res}\limits_{\lambda =\frac{n}{2} } R_{11} \left(x,t,\lambda \right)} \end{array}\]
Separately calculate every item.

\[R_{11} \left(x,t,\lambda \right)-R_{12} \left(x,t,\lambda \right)=\frac{f_{1}^{+} \left(x,\lambda \right)f_{1}^{+} \left(t,\lambda \right)}{2i\lambda C_{12} \left(\lambda \right)C_{22} \left(\lambda \right)} \]

\noindent  Residues of resolvent $R_{11} \left(x,t,\lambda
\right)\, \, \, $ in $\lambda _{1,} \lambda _{2,} ....\lambda
_{l}\,\,\, $denote by $G_{11} \left(\lambda _{n} \right)$ . Then
$G_{11} \left(\lambda _{n} \right)$ will be equal to

\[G_{11} \left(\lambda _{n} \right)=\mathop{\lim }\limits_{\lambda \to \lambda _{n} } \left(\lambda -\lambda _{n} \right)R_{11} \left(x,t,\lambda \right). \]
Then for every function $\psi \left(x\right)$ belonging to $L_{2} \left(-\infty ,+\infty ,\rho \left(x\right)\right)$ we get following eigenfunction expansion in the form

\begin{equation} \label{GrindEQ__25_}
\begin{array}{l} {\psi \left(x\right)=-\frac{1}{2i\pi } \int _{-\infty }^{+\infty }\rho \left(t\right)\psi \left(t\right) [\oint _{\Gamma_{0}^{-} }\frac{f_{1}^{+} \left(x,\lambda \right)f_{1}^{+} \left(t,\lambda \right)}{2i\lambda C_{12} \left(\lambda \right)C_{22} \left(\lambda \right)}  ]d\lambda +} \\ {+G_{11} \left(\lambda _{n} \right)+\frac{2}{in} V_{nn} f_{1}^{+} \left(x,\frac{n}{2} \right)f_{1}^{+} \left(t,\frac{n}{2} \right)+F\left(x,t\right)]dt} \end{array}
\end{equation}

\noindent

\noindent \textbf{\underline{SOLUTION OF THE INVERSE PROBLEM}.}

\noindent \textbf{}

Let's study the inverse problem for the problem (1-3). In spectral
expansion \eqref{GrindEQ__25_} the numbers \textbf{$V_{nn} $ }play
a part of normalizing numbers for the function $f_{1}^{+}
\left(x,\lambda \right)$ responding it spectral singularities.

The inverse problem is formulae as follows.

\noindent \textbf{\underline{INVERSE PROBLEM}. }Given the spectral
data $\{ \, \, \, C_{12} \left(\lambda \right),V_{nn} \} $
construct the $\beta $  and potential $q\left(x\right)$.

\noindent Using the results obtained above we arrive at the following procedure for solution of the inverse problem.

\noindent 1. Taking into account \eqref{GrindEQ__6_} we get

\[V_{n,\alpha +n} =V_{nn} \sum _{m=1}^{\alpha }\frac{V_{m\alpha } }{m+n}  ,\]
from which all numbers  $V_{n\alpha } ,\, \, \, \alpha =1,2.....,\, \, \, n=1,2,....n<\alpha $ are defined.

\noindent 3. Then from recurrent formula
\eqref{GrindEQ__6_}-\eqref{GrindEQ__8_}, find all numbers
$q_{n}$.

\noindent 4. The number $\beta $  is defined from equality

\[\beta =-2\mathop{\lim }\limits_{Im\lambda \to \infty } C_{12} \left(\lambda \right)-1.\]

So inverse problem has a unique solution and the numbers $\beta $ and $q_{n} $ are defined constructively by the spectral data.

\noindent \textbf{\underline{Theorem 2}.} The specification of the
spectral data uniquely determines  $\beta $ and potential
$q\left(x\right)$.

\noindent \textbf{\underline{REFERENCES}.}

\noindent 1.   Kamke E. Handbook of ordinary Differential
equations (Russian).  Nauka,

\noindent       Moscow,1976.

\noindent 2.   Jeffrey C. Lagarias. The Schrodinger Operator with
Morse Potential on the

\noindent       right half line. Arxiv: 0712.3238v1 [math.SP] 19 Dec 2007.

\noindent 3.   P. C. Sabatier and B. Dolveck-Guilpard.1998 On
modelling discontinuous

\noindent        media. One-dimensional approximations, J.Math. Phys. 29, 861-868.

\noindent 4.   P. C. Sabatier.1989 On modelling discontinuous
media. Three-dimensional

\noindent         scattering, J. Math. Phys. 30, 2585-2598 .

\noindent 5.   F. Dupuy and P C. Sabatier.1992 Discontinuous media and undetermined

\noindent         scattering problems, J. Phys. A 25, 4253-4268.

\noindent 6.    F R Molino and P C. Sabatier.1994 Elastic waves in discontinuous media:

\noindent        Three-dimensional scattering,J. Math. Phys.3, 45944635.

\noindent 7.    T. Aktosun, M. Klaus, and C. van der Mee,1993 On the Riemann-Hilbert

\noindent         problem for the one-dimensional Schrodinger
                  equation, J. Math. Phys. 34,2651-2690

 \noindent8.   Guseinov, I.~M. and Pashaev, R.~T.~ On an inverse
               problem for a second-order differential equation. UMN, 2002 , 57:3, 147--148

\noindent 9.   Gasymov M.G. Spectral analysis of a class non-self-adjoint operator of the

\noindent       second order. Functional analysis and its appendix. (In  Russian) 1980,

\noindent       V34.1.pp.14-19

\noindent 10. Pastur L.A., Tkachenko V.A. An inverse problem for one class of

\noindent       onedimentional Shchrodinger's operators with complex periodic potentials.

\noindent       Funksional analysis and its appendix.(in Russian)
1990 V54.¹6.pp. 1252-1269

\noindent11.  V.Guillemin, A. Uribe. Hardy functions and the
               inverse spectral method.

\noindent       Comm.In Partial Differential equations, 8\eqref{GrindEQ__13_},1455-1474(1983)

\noindent 12.   Efendiev, R.F. Spectral analysis of a class of non-self-adjoint differential

\noindent       operator pencils with a generalized function.Teoreticheskaya i

\noindent       Matematicheskaya  Fizika, 2005,Vol.145,1.pp.102-107,October(Russian).

\noindent       Theoretical and Mathematical Physics,145\eqref{GrindEQ__1_}:1457-461,(English).

\noindent 13.  Efendiev, R. F. Complete solution of an inverse problem for one class of the

\noindent       high order ordinary differential operators with periodic coefficients. Zh. Mat.

\noindent       Fiz. Anal. Geom2006, C.\textit{ }2, no. 1, 73--86, 111.

\noindent 14. Efendiev, R.F. The Characterization Problem for One Class of Second Order

\noindent       Operator Pencil with Complex Periodic Coefficients. Moscow Mathematical

\noindent       Journal, 2007, Volume~7,~Number~1, pp.55-65.

\end{document}